\documentstyle[12pt]{article}
\def\c#1#2#3{c_{#1#2}{}^{#3}}

\def\r#1{(\ref{#1})}

\def\bel{\begin{equation}\label}
\def\ee{\end{equation}}
\def\ba{\begin{array}}
\def\ea{\end{array}}

\def\to{\rightarrow}

\def\cop{\triangle}

\def\a{\alpha}
\def\grad{\mathop{\rm grad}}

\begin{document}
\def\tens{\otimes}
\def\ben{\begin{enumerate}}
\def\een{\end{enumerate}}
\def\g#1#2{\grad(#1,#2)}
\def\hg{\hat g}
\def\<{\left<}
\def\>{\right>}
\title{
\rightline{{\raisebox{1cm}{\normalsize \bf IFT UWR 918/97}}}
Classification of low dimensional Lie super-bialgebras
\author{C. Juszczak $^*$ and J. T. Sobczyk
\thanks{Supported by KBN grant 2P 30208706}}
}
\date{\normalsize
Institute of Theoretical Physics, University of Wroc\l{}aw,\\
\normalsize Pl. Maksa Borna 9,
Wroc\l{}aw, Poland\vspace*{-3mm}}
\maketitle
\begin{abstract}
A thorough analysis of Lie super-bialgebra structures
on Lie
super-algebras $osp(1;2)$ and super $e(2)$ is presented.
Combined technique of computer algebraic computations and a subsequent
identification of equivalent structures is applied. In all the cases
Poisson-Lie brackets on supergroups are found. Possibility of
quantizing them in order to obtain quantum groups is discussed. It turns
out to be straightforward for all but one structures for super-$E(2)$
group.
\end{abstract}

\section{Introduction}\label{sec1}

About 10 years of intensive research in quantum groups
\cite{FRT}
\cite{D1}
demonstrated
richness of investigated mathematical structures. No wonder that except
from efforts to find their physical applications many authors tried
to classify them from a mathematical point of view. The classifications are
based on many different ideas. We do not wish to provide exhaustive list
of them but would like to mention some of well known examples.
If quantum $SL(2)$ group is defined as one
that preserves a non-degenerate bilinear form two inequivalent deformations
appear
\cite{D-V}.
More complicated analysis of the group $GL(3)$ leads to a conclusion that $26$
inequivalent deformations exist in this case
\cite{EO}.
Ref. \cite{Ohn} contains a very recent attempt to classify quantum $SL(3)$
groups.
A classification of all the $4\times 4$ matrices satisfying
Yang-Baxter equation  was done in
\cite{Hiet}.
Later this result was used
in order to list possible quantum deformations of the supergroup
$GL(1,1)$
\cite{FHR}.
An interest in deformations of $D=4$ relativistic
symmetries led to a classification of deformations of Lorentz
and Poincar\' e groups
\cite{WPZ}.

{}From the early years of interest in quantum groups it is
well known that
they are closely related with notions of Lie bialgebra and Lie-Poisson
group. Suppose that
a parameter $q$ is introduced such that the value $q=1$ corresponds
to undeformed universal enveloping algebra $U({\cal G})$ or a classical
commutative
algebra of functions on a Lie group $G$ (with the Lie algebra ${\cal G}$). Then
the structure of Hopf algebra of the quantum group gives rise in the
classical limit to extra structure on $G$ or ${\cal G}$. $G$ becomes a
Poisson-Lie group and ${\cal G}$ a Lie bialgebra
\cite{D2}.
These structures
can be viewed as possible {\it directions}
of quantum deformations. A fundamental
fact about them is that all can be quantized i.e. they arise as
the classical limit from a bona fide quantum deformation
\cite{EK}
This
suggests that an effective approach to classify quantum groups could be by
classifying their classical limits and that seems to be much easier to
be pulled off. The
study of Lie bialgebras becomes in some cases even easier as they might be
of coboundary type what happens if they can be described my means of a classical
$r$-matrix satisfying (modified) classical Yang-Baxter equation. Cohomological
properties of semisimple Lie groups imply that all related
Lie bialgebras are in fact coboundaries
\cite{D2}.
More elaborate argument shows that the same
is true in the case of inhomogeneous groups of space time symmetries for
any signature of metric for dimensionality of space-time $D>2$
\cite{Zak1}.
Thus the problem of
classification of Lie bialgebra structures contains as a subproblem
a classification
of solutions of (modified) classical Yang-Baxter
equation which has been studied by many authors.
It should be stressed
however
that in most cases Lie algebras admit bialgebra structures which are
not coboundaries.
Recently many low
dimensional Lie algebras were investigated for bialgebra structures
revealing a surprising number of possibilities. Two-dimensional Galilei
algebra admits 9 inequivalent Lie bialgebra structures and
only one out of them is a coboundary
\cite{Kowal}.
The central extension of Galilei algebra by the
mass operator admits 26 inequivalent Lie bialgebra structures out of which 8
are coboundaries
\cite{Opan}.
Similar analysis aimed mainly on the construction of quantum groups was
performed before in the cases of Heisenberg-Weyl and oscillator Lie
algebras
\cite{Hiszp}.
The classification of the possible classical $r$-matrices (and
automatically of Lie bialgebra structures) was done in the
case of Lorentz and Poincar\' e algebras
\cite{Zak2}.
Classical $r$-matrices for $SL(3)$ were listed in \cite{Stolin}. Corresponding
$R$-matrices satisfying Yang-Baxter equation were found in \cite{GG}.

In this paper we will study quantum deformations of supergroups and their
classical limits - Lie super-bialgebras. Quantum supergroups have been
studied by many authors.
Knowledge of $R$-matrices and low-dimensional representations
of $osp(1|2)$, $su(1|1)$, $gl(2|1)$, $sl(1|2)$
has been
used to construct integrable models
\cite{Zastos}.
Up to our knowledge no systematic investigation of Lie super-bialgebras
has been yet undertaken. We decided to study two cases: $osp(1,2)$ and
supersymmetric extension of $e(2)$ algebra both treated as complex.
The first one is interesting as
it plays in the supersymmetric case
a role analogous to $sl(2)$. Many papers devoted
to quantum supergroups are based on a quantum deformation proposed by
\cite{Kulish}.
A natural question is: are other deformations possible like it is
in the case of $sl(2)$?
The
supersymmetric extension of $e(2)$ was chosen as $e(2)$ group is one
which has a simple structure but admits several inequivalent Lie
bialgebras
\cite{Zak3}
\cite{Sobczyk}.
There are 4 of them  with one being strictly
speaking a 1-parameter family. Two structures are generated by classical
$r$-matrices. All the four quantum deformations of the group $E(2)$ are
known which could be helpful in finding quantum deformations in the
supersymmetric case.
Some quantum deformations of super-$e(2)$ algebra have been already
discussed in the literature. They can be (in analogy to $e(2)$ case)
obtained by means of a contraction procedure
\cite{Kontrakcja}.

The problem of classifying Lie bialgebra structures for a given Lie algebra
is mostly a computational one. We found it effective to use at various
steps a computer. For a Lie algebra with $N$ generators it is necessary to
deal with equations containing the number of parameters which increases
as $N^3$. Some of the equations are linear but remaining (co-Jacobi
identities) are quadratic. The complexity of the problem grows quickly
with $N$. The biggest obstacle in continuing with this program for large
$N$ is primarily the fact
that not everything can be done by computer. The most delicate
part of the problem is to single out orbits under the action of the group of
automorphisms of ${\cal G}$ in a computer produced space of solutions.

A part of our results (a list of bialgebra structures but not of Poisson-Lie
brackets) was presented in \cite{Juszczak}.
The first result is that in the $osp(1,2)$ case all the Lie super-bialgebras
are coboundaries. We found three
independent classical $r$-matrices. One solution
is quite obvious. As $sl(2)$ is a subalgebra of $osp(1,2)$ any $r$-matrix for
the former algebra satisfying (not modified!) classical YB equation is
automatically a $r$-matrix for the latter. The second solution corresponds
to the quantum deformation of $osp(1,2)$ described in
\cite{Kulish}. Relations for quantum $OSp(1,2)$ are given in \cite{KCh}.
The third
solution is a new one and requires a detailed study.

In the case of super-$e(2)$ we found
six families of independent Lie super-bialgebra
structures. In their most general form they are not coboundaries  but four
of them contain coboundary members. Some links to the
classification of analogous structures on $e(2)$ can be established.
The case (i) in our list
is in a clear analogy with $\delta_2$ of
\cite{Sobczyk}.
The case (iv)
corresponds to $\delta_1$.
It is the infinitesimal form of quantum deformation introduced in
\cite{Kontrakcja}.
The existence of $r$-matrix $r_2$ is again obvious
as it consists only of generators of $e(2)$ and satisfies (not modified!)
classical YB equation.

For all the Lie super-bialgebra structures we calculated corresponding
Lie-Poisson brackets. Their form is such that in five cases of
super $E(2)$ group it is easy
to go through with the program of constructing quantum supergroups
\cite{Zak4}
\cite{KLMS}.
It is sufficient to change Lie-Poisson brackets into (anti) commutators.

The paper is organized as follows. In Chapter \ref{sec2} all the basic
concepts and notation used in the rest of the paper are introduced.
Chapters \ref{sec3} and \ref{sec4} contain derivation
of inequivalent Lie super-bialgebras for $osp(1,2)$ and
super-$e(2)$. In Chapter \ref{sec5}
 Poisson-Lie brackets for both cases are presented. In Chapter \ref{sec6}
  our conclusions and some final remarks are presented.

After completing our study we have learned about the paper 
in which by applying
Drinfeld twisting procedure new
quantum supergroup structure on $OSp(1,2)$ was found \cite{CeKu}.
It is clear that it
corresponds to the above mentioned
"trivial" Lie super-bialgebra structure (number 1 on our
list presented in Chapter \ref{sec4}).
\section{Basic definitions and notation}\label{sec2}
Super Lie-algebra ${\cal G}$ is a graded vector space
\bel{ca}
{\cal G}={\cal G}_B \oplus {\cal G}_F
\ee
with the grade function $grade({\cal G}_B) =0$, $grade({\cal G}_F)=1$.
A Lie superalgebra structure is provided by a
linear mapping
\bel{cb}
[\ ,\ ] : {\cal G} \tens {\cal G} \to {\cal G}
\ee
satisfying requirements of (graded) antisymmetry and Jacobi identity.
In order to express them it is useful to introduce a basis in ${\cal G}$
$\{g_i\} \subset {\cal G}_B \cup {\cal G}_F$ and structure constants
\bel{cc}
[g_i,g_j]= c_{ij}^{\ \ k} g_k.
\ee
Structure constants have to satisfy
\bel{cca}
c_{ij}^{\ \ k}=0\ \  \mbox{ whenever }\ \  grade (g_i) + grade (g_j) \neq
grade (g_k) \mbox{ (mod 2)}
\ee
\bel{ccaa}
c_{ij}^{\ \ k}= - z(i,j) c_{ji}^{\ \ k},
\ee
and
\bel{ccb}
c_{ij}^{\ \ k} c_{kl}^{\ \ m} z(i,l)
+ c_{jl}^{\ \ k} c_{ki}^{\ \ m} z(j,i)
+ c_{li}^{\ \ k} c_{kj}^{\ \ m} z(l,j) =0
\ee
where
\bel{ce}
z(i,j) = (-1)^{grad (g_i)grad (g_j)}.
\ee

\noindent
Lie super-bialgebra structure is a linear mapping
\bel{cf}
\delta : {\cal G} \to {\cal G} \tens {\cal G}
\ee
which in the chosen basis reads
\bel{cg1}
\delta( g_i)= f_i^{\ kl} g_k \tens g_l.
\ee
$\delta$ has to satisfy several requirements. First of all it makes
the dual linear space ${\cal G}^*$ a Lie super-algebra
\bel{cea}
f_k^{\ ij} = 0\ \  \mbox{ whenever }\ \ grad (g_i) + grad (g_j) \neq
grad (g_k)\mbox{ (mod 2)}
\ee

\bel{ceab}
f_k^{\ ij} = - z(i,j) f_k^{\ ji}
\ee

\bel{ceac}
f_i^{\ kj} f_j^{\ lm} z(k,m) + f_i^{\ lj} f_j^{\ mk} z(l,k) +
f_i^{\ mj} f_j^{\ kl} z(m,l) =0.
\ee
Moreover structure constants $c$'s and $f$'s have to be related
\bel{cg}
c_{ij}^{\ \ k} f_k^{\ lm} = f_i^{\ lk} c_{kj}^{\ \ m}
+ c_{kj}^{\ \ l} f_i^{\ km} z(m,j)
+ c_{jk}^{\ \ l} f_j^{\ km} + f_j^{\ lk} c_{ik}^{\ \ m} z(i,l)
\ee

\noindent
Coboundary Lie super-bialgebra is a pair $({\cal G},r)$, where ${\cal G}$
is a Lie
super-bialgebra and $r \in {\cal G}_B \wedge {\cal G}_B
\oplus {\cal G}_F \wedge {\cal G}_F \subset {\cal G}
\wedge {\cal G}$.
such that for every $g_i \in {\cal G}$
\bel{ch}
\delta(g_i) = [r, g_i\tens 1 + 1 \tens g_i]
\ee

\noindent
Schouten bracket is defined as follows
\bel{5}
[[r,r]] \equiv [r_{12},r_{13}] + [r_{12},r_{23}] + [r_{13},r_{23}]
\ee
where $r_{12} = r \tens 1$, $r_{23} = 1 \tens r$, ....

\noindent
$r$ satisfies classical Yang-Baxter equation (CYBE) if
\bel{ci}
[[r,r]] = 0
\ee
and modified CYBE  if
\bel{cj}
\forall {g_i\in {\cal G}}\quad [[[r,r]], g_i\tens1\tens1 + 1\tens g_i
\tens1+1\tens1\tens g_i]=0.
\ee

\noindent
Super $e(2)$ Lie algebra is spanned by the set of generators
$\{ H$, $P_+$, $P_-$, $D_+$, $D_-  \}$ which fulfill the following (anti)
commutation relations
\bel{8}
\ba{rcl}
[H, P_\pm] &=& \pm P_\pm\,,\\{}
[H, D_\pm] &=& \pm \frac12 D_\pm\,,\\{}
[P_+, P_-] &=& 0\,,\\{}
\{D_+,D_-\} &=& 0\,,\\{}
\{D_\pm,D_\pm\} &=& P_\pm\,,\\{}
[P_\pm, D_\pm] &=& 0\,,\\{}
[P_\pm, D_\mp] &=& 0\,,
\ea
\ee

A convenient parameterization of the classical super-$E(2)$ group is
obtained by means of exponentiation
\bel{aa}
g(s,a,b,\xi ,\eta ) = \exp (sH) \exp (aP_+)\exp (bP_-)\exp (\xi D_+)
\exp (\eta D_-).
\ee
Coproducts are
\bel{ab}
\Delta (s) = s\otimes 1 + 1\otimes s,
\ee
\bel{ac}
\Delta (a) = 1\otimes a + a\otimes \exp (-s) + {1\over 2} \xi\otimes
\xi\exp (-{s\over 2}),
\ee
\bel{ad}
\Delta (b) = 1\otimes b + b\otimes \exp (s) + {1\over 2} \eta\otimes
\eta\exp ({s\over 2}),
\ee
\bel{ae}
\Delta (\xi ) = 1\otimes \xi + \xi\otimes \exp (-{s\over 2}),
\ee
\bel{af}
\Delta (\eta ) = 1\otimes \eta + \eta\otimes \exp ({s\over 2}).
\ee

\noindent
Lie superalgebra $osp(1,2)$ is spanned by the set of generators
$\{ H$, $X_+$, $X_-$, $V_+$, $V_-  \}$ which fulfill the following (anti) commutation
relations
\bel{8a}
\ba{rcl}
[H, X_\pm] &=& \pm X_\pm\,,\\{}
[H, V_\pm] &=& \pm \frac12 V_\pm\,,\\{}
[X_+, X_-] &=& 2H\,,\\{}
\{V_+,V_-\} &=& -\frac12 H\,,\\{}
\{V_\pm,V_\pm\} &=& \pm\frac12 X_\pm\,,\\{}
[X_\pm, V_\pm] &=& 0\,,\\{}
[X_\pm, V_\mp] &=& V_\pm\,,
\ea
\ee
In the case of supergroup $OSp(1|2)$ we will use more implicit
parameterization by means of $3\times 3$ supermatrices subject to
certain constraints
\cite{Kulish}

\bel{ag}
g = \pmatrix{a&\alpha&b\cr \gamma&e&\beta\cr c&\delta&d}
\ee
where
$e=1+\alpha\delta$,\ \ \  $\gamma = c\alpha - a\delta$, \ \ \
$\beta = d\alpha - b\delta$, \ \ \ $ad - bc + \alpha\delta =1$.
Variables denoted by Greek letters are of Grassmanian type.
Coproducts follow from matrix multiplication of elements of $G$.
\section{Lie super-bialgebras for $osp(1,2)$}\label{sec3}

The problem studied in Chapters \ref{sec3} and \ref{sec4}
can be formulated
in the following way:

\medskip
{\sl Given a set of structure constants $c_{ij}^{\ \ k}$
find all sets of
structure constants $f_m^{\ np}$
that give rise to Lie super-bialgebras.
Two such sets
are considered equivalent if they can be made equal by a
change of the basis in ${\cal G}$.

Verify if obtained structures
are or not of coboundary type.}

\medskip
Initial steps in the analysis can be made using a computer.
An arbitrary form of $f$'s satisfying super-antisymmetry
and preserving the grading is assumed. Constraints from
the set of linear equations coming from the cocycle
condition \r{cg} are firstly taken into account. Then
the set of quadratic equations coming from the
super co-Jacobi conditions \r{ceac} is to be solved.

In the case of $osp(1,2)$
we discover that all the possible bialgebra structures
are coboundaries and that the classical $r$-matrix is in one of two
possible forms:
\bel{9}
\ba{rcl}
r_a &=& x(X_+ \wedge X_- + 2 V_+ \wedge V_-) \\
    &&+ y(H \wedge X_+ -  V_+ \wedge V_+) \\
    &&+ z(H \wedge X_- - V_- \wedge V_-)\,, \\[3mm]
r_b &=& \pm \sqrt{uv} X_+ \wedge X_- + u H \wedge X_+ + v H \wedge X_- \,,
\ea
\ee
where $x$, $y$, $z$, $u$, $v$ are arbitrary complex numbers.

\noindent
Identification of automorphisms of superalgebra $osp(1,2)$ is
fairly straightforward.
They can not mix fermions with bosons
(grading is preserved). As the superalgebra is generated just by the two
fermions $V_+$ and $V_-$ every automorphism
is determined by the
following transformation of the fermions;
\bel{10}
\ba{rcl}
\widetilde V_+ &=& a\, V_+ + b\, V_-\\
\widetilde V_- &=& c\, V_+ + d\, V_-.
\ea
\ee
We easily derive that under the above transformation
\bel{11}
\ba{lcl}
\widetilde H &=& -ac\, X_+ + (ad +bc )\, H + bd\, X_-\,,\\
\widetilde X_+ &=& \,\,\,\,a^2\, X_+ -2 ab\, H - b^2\, X_-\,,\\
\widetilde X_- &=& -c^2\, X_+ + 2cd\,  H + d^2\, X_-\,.
\ea
\ee
The operators $\{\widetilde H,\widetilde X_+, \widetilde X_-, \widetilde
V_+,\widetilde V_-\}$ obey the same super-commutation relations
as $\{ H, X_+,  X_-, 
V_+, V_-\}$
if and only if
\bel{12}
\det \pmatrix {a&b \cr c &d } = 1\,.
\ee


\noindent
Under the transformation \r{10}-\r{11} the parameters $x$, $y$, $z$ 
of the $r$-matrix $r_a$ given in \r{9} transform as follows
\bel{13a}
\pmatrix{ \widetilde y & -\widetilde x \cr
-\widetilde x & \widetilde z } = 
\pmatrix{ a & b \cr
 c & d }  
\pmatrix{ y & -x \cr
-x & z }  
\pmatrix{ a & c \cr
 b & d }\,.
\ee
We notice that this is exactly the way the symmetric form transforms under 
the change of basis. If we take into account the Sylvester theorem we see
that we can make the matrix 
$$\pmatrix{ y & -x \cr -x & z} $$
diagonal ($x=0$) with $y=z$ or $y=1$, $z=0$ (without the condition \r{12} it would 
be possible to make either $y=z=1$ or $y=1$, $z=0$).

Thus $r_a$ if different from zero is equivalent to one of the 
following:
\bel{14}
r_a = \cases{r_2\equiv H \wedge X_+ - V_+ \wedge V_+ &  if $ x^2-yz =0$ \cr
       r_3 \equiv   t(H \wedge X_+ - V_+ \wedge V_+ +H \wedge X_- - V_- \wedge 
V_-) & if $x^2-yz \neq 0$\,.}
\ee
By means of analogous reasoning we arrive at the conclusion that the
$r$-matrix $r_b$ is always equivalent to
\bel{15}
r_1\equiv H\wedge X_+\,.
\ee

\section{Lie super-bialgebras for  super-$e(2)$}\label{sec4}

\noindent
The initial steps of the analysis are made using a computer like in the case
of $osp(1,2)$. After solving equations quadratic in structure constants it
turns out that most of bialgebra structures are not coboundaries.
This makes the investigation more involved since we cannot
use the $r$-matrix formulation and must explicitly state the co-Lie structure
for each solution found.

\noindent
The possibilities found by means of the computer are:
\ben
\item[Case \bf A] 
\bel{16}
\ba{lcl}
\delta(H) &=& H \wedge (aP_+ + bP_-) + c P_+ \wedge P_-\,,\\[1mm]
\delta(P_+)&=& P_+ \wedge b P_-\,,\\[1mm]
\delta(P_-)&=& -a P_+ \wedge P_-\,,\\[1mm]
\delta(D_+)&=& \frac12 (a P_+ -b P_-) \wedge D_+ \pm \sqrt{ab} P_+ \wedge 
D_- \,,\\[1mm]
\delta(D_-)&=& \frac12 (a P_+ -b P_-) \wedge D_- \pm \sqrt{ab} P_- \wedge 
D_+ \,,
\ea
\ee
\item[Case \bf B]
\bel{17}
\ba{lcl}
\delta(H) &=&~ a( H \wedge P_+  -\frac12 D_+ \wedge D_+) \\
&&+ b(H \wedge P_- +\frac12 D_-\wedge D_-) +c P_+ \wedge
P_-\,,\\[1mm]
\delta(P_+)&=& P_+ \wedge b P_- +d(2H\wedge P_+ - D_+\wedge D_+)\,,\\[1mm]
\delta(P_-)&=& -a P_+ \wedge P_- +d(2H\wedge P_- + D_-\wedge D_-)\,,\\[1mm]
\delta(D_+)&=& -\frac12 (a P_+ + b P_-) \wedge D_+ +d( H \wedge D_+ )\,,\\[1mm]
\delta(D_-)&=& -\frac12 (a P_+ + b P_-) \wedge D_-+ +d( H \wedge D_-)\,,
\ea
\ee
where $a$, $b$, $c$, $d$ are arbitrary complex numbers such that $cd=0$.
\een

The set of automorphisms of super-$e(2)$ is generated by three transformations:
\ben
\item $\widetilde H = H + \alpha P_+ + \beta P_-\,,\qquad
\widetilde P_\pm = P_\pm\,,\qquad \widetilde D_\pm = D_\pm\,,$ 

\item $\widetilde H = -H\,,\qquad 
\widetilde P_\pm = P_\mp\,,\qquad \widetilde D_\pm = D_\mp\,,$ 

\item $\widetilde H = H\,,$ \\
$\widetilde P_+ = \alpha^2 P_+\,,\qquad \widetilde D_+ =\alpha 
D_+\,,$ \\
$\widetilde P_- = \beta^2 P_-\,,\qquad \widetilde D_- =\beta
D_-\,,$ 
\een

\noindent
When taken into account they lead to a conclusion
that $a$ and $b$ can be
scaled out to take value of $1$ or $0$. It is also possible to show
that when $d\neq0$ then by taking $\widetilde H = H + \frac{a}{2d} P_+ + 
\frac{b}{2d} P_-$, we can make $a=b=c=0$.

\noindent
Finally we arrive at six families of bialgebra structures
whose representatives can be described by means of the following substitutions:
\ben
\item[(i)] Case {\bf A} where $a=b=0$,
\item[(ii)] Case {\bf A} where $a=1$, $b=0$,
\item[(iii)] Case {\bf A} where $a=b=1$,
\item[(iv)] Case {\bf B} where $a=b=c=0$,
\item[(v)] Case {\bf B} where $a=1$, $b=d=0$,
\item[(vi)] Case {\bf B} where $a=b=1$, $d=0$.
\een

\noindent
Classical $r$-matrices exist if both $c$ and $d$ vanish. 
Their general forms for cases {\bf A} and {\bf B} are
\bel{19aa}
\ba{rcl}
r_A &=& a H\wedge P_+ - b H \wedge P_- + \sqrt{ab} D_+ \wedge D_- + f P_+
\wedge P_-\,,\\
r_B &=& a (H\wedge P_+ -\frac12 D_+ \wedge D_+) - b (H \wedge P_- + \frac12
D_-\wedge D_-) + f P_+ \wedge P_-\,,
\ea
\ee
Term $P_+\wedge P_-$ is irrelevant and will be omitted
since its commutators \r{ch} with all the
generators of super-$e(2)$ vanish.

\noindent
After substitutions we obtain
\bel{19}
\ba{rcl}
r_{(ii)} &=&  H\wedge P_+\,,\\
r_{(iii)} &=&  H\wedge P_+ -  H \wedge P_- +  D_+ \wedge D_-  \,,\\
r_{(v)} &=&  H\wedge P_+ -\frac12 D_+ \wedge D_+\,,\\
r_{(vi)} &=&  H\wedge P_+ -\frac12 D_+ \wedge D_+ - H \wedge P_- - \frac12
D_-\wedge D_- \,,
\ea
\ee
Only $r_{(ii)}$ and $r_{(v)}$ satisfy CYBE.

\section{Poisson-Lie brackets}\label{sec5}

\noindent
Lie super-bialgebras are in $1:1$ correspondence with Poisson-Lie
structures on supergroups
\cite{Andrusz}. 
Poisson-Lie brackets satisfy the following properties:
\begin{eqnarray}
\label{-a1111}
\{f,g\}&=& -z(f,g)\{g,f\}\,,
\\
\label{a2222}
\{f,gh\}&=& \{f,g\}h+z(f,g)g\{f,h\}\,,
\\
\label{a3333}
0&=&z(f,h)\{f,\{g,h\}\} +
z(g,f)\{g,\{h,f\}\}\nonumber \\ 
&&+ z(h,g)\{h,\{f,g\}\}\,,\\ 
\label{a4444}
\cop\{f,g\} &=& \{\cop f,\cop g\}\,.
\end{eqnarray}

Poisson-Lie brackets are most easily introduced by means of
left- and right-invariant vector fields on a supergroup. In
the super case one should distinguish left- and right-hand side derivatives
of superfunctions. In coboundary case with a classical $r$-matrix $r_{kl}$
Poisson-Lie brackets are given by
\cite{KLMS}

\bel{ba}
\{\phi , \psi\} = \left( Y_k^{(r)}\phi\right)
r_{kj} \left( Y_j^{(l)}\psi\right) -
\left( X_k^{(r)}\phi\right) r_{kj} \left( X_j^{(l)}\psi\right)
\ee

\noindent
where $Y_k^{(r,l)}$ denotes left-invariant right (r) or left (l) derivatives
and $X_k^{(r,l)}$ right-invariant right (r) or left (l) derivatives. $Y_k$
and $X_k$ can be derived from the coproducts.

We present below how they
act on generators of both supergroups.

\ben
\item[a)] super-$E(2)$
\bel{bb}
Y_H^{(r,l)}\pmatrix{a\cr b\cr s\cr \xi\cr\eta} =
\pmatrix{-a\cr b\cr 1\cr -\xi /2 \cr \eta /2},\qquad
X_H^{(r,l)}\pmatrix{a\cr b\cr s\cr \xi\cr\eta} =
\pmatrix{0\cr 0\cr 1\cr 0 \cr 0},
\ee

\bel{bc}
Y_{P_+}^{(r,l)}\pmatrix{a\cr b\cr s\cr \xi\cr\eta} =
\pmatrix{1\cr 0\cr 0\cr 0 \cr 0},\qquad
X_{P_+}^{(r,l)}\pmatrix{a\cr b\cr s\cr \xi\cr\eta} =
\pmatrix{\exp (-s)\cr 0\cr 0\cr 0 \cr 0},
\ee

\bel{bd}
Y_{P_-}^{(r,l)}\pmatrix{a\cr b\cr s\cr \xi\cr\eta} =
\pmatrix{0\cr 1\cr 0\cr 0 \cr 0},\qquad
X_{P_-}^{(r,l)}\pmatrix{a\cr b\cr s\cr \xi\cr\eta} =
\pmatrix{0\cr \exp {(s)} \cr 0\cr 0 \cr 0},
\ee

\bel{be}
Y_{D_-}^{(r)}\pmatrix{a\cr b\cr s\cr \xi\cr\eta} =
\pmatrix{0\cr \eta /2\cr 0\cr 0 \cr 1},\qquad
X_{D_-}^{(r)}\pmatrix{a\cr b\cr s\cr \xi\cr\eta} =
\pmatrix{0\cr -\eta \exp ({s\over 2}) /2\cr 0\cr 0 \cr \exp ({s\over 2})},
\ee

\bel{bf}
Y_{D_-}^{(l)}\pmatrix{a\cr b\cr s\cr \xi\cr\eta} =
\pmatrix{0\cr -\eta /2\cr 0\cr 0 \cr 1},\qquad
X_{D_-}^{(l)}\pmatrix{a\cr b\cr s\cr \xi\cr\eta} =
\pmatrix{0\cr \eta \exp ({s\over 2}) /2\cr 0\cr 0 \cr \exp ({s\over 2})},
\ee

\bel{bg}
Y_{D_+}^{(r)}\pmatrix{a\cr b\cr s\cr \xi\cr\eta} =
\pmatrix{\xi /2\cr 0 \cr 0\cr 1 \cr 0},\qquad
X_{D_+}^{(r)}\pmatrix{a\cr b\cr s\cr \xi\cr\eta} =
\pmatrix{-\xi\exp (-{s\over 2}) /2\cr 0\cr 0\cr \exp (-{s\over 2}) \cr 0},
\ee

\bel{bh}
Y_{D_+}^{(l)}\pmatrix{a\cr b\cr s\cr \xi\cr\eta} =
\pmatrix{-\xi /2\cr 0 \cr 0\cr 1 \cr 0},\qquad
X_{D_+}^{(l)}\pmatrix{a\cr b\cr s\cr \xi\cr\eta} =
\pmatrix{\xi\exp (-{s\over 2}) /2\cr 0\cr 0\cr \exp (-{s\over 2}) \cr 0}.
\ee

\item[b)] $OSp(1|2)$

\bel{bi}
Y_H^{(r,l)}\pmatrix{a\cr \alpha\cr b\cr c\cr \delta\cr d} =
\pmatrix{a/2\cr 0\cr -b/2\cr c/2\cr 0\cr -d/2},\qquad
X_H^{(r,l)}\pmatrix{a\cr \alpha\cr b\cr c\cr \delta\cr d} =
\pmatrix{a/2\cr \alpha /2\cr b/2\cr -c/2\cr -\delta /2\cr -d/2},
\ee

\bel{bj}
Y_{X_+}^{(r,l)}\pmatrix{a\cr \alpha\cr b\cr c\cr \delta\cr d} =
\pmatrix{0\cr 0\cr a\cr 0\cr 0\cr c},\qquad
X_{X_+}^{(r,l)}\pmatrix{a\cr \alpha\cr b\cr c\cr \delta\cr d} =
\pmatrix{c\cr \delta\cr d\cr 0\cr 0\cr 0},
\ee

\bel{bk}
Y_{X_-}^{(r,l)}\pmatrix{a\cr \alpha\cr b\cr c\cr \delta\cr d} =
\pmatrix{b\cr 0\cr 0\cr d\cr 0\cr 0},\qquad
X_{X_-}^{(r,l)}\pmatrix{a\cr \alpha\cr b\cr c\cr \delta\cr d} =
\pmatrix{0\cr 0\cr 0\cr a\cr \alpha\cr b},
\ee

\bel{bl}
Y_{V_+}^{(r)}\pmatrix{a\cr \alpha\cr b\cr c\cr \delta\cr d} =
\pmatrix{0\cr a/2\cr \alpha /2\cr 0\cr c/2\cr \delta /2},\qquad
X_{V_+}^{(r)}\pmatrix{a\cr \alpha\cr b\cr c\cr \delta\cr d} =
\pmatrix{-\gamma /2\cr e/2\cr -\beta /2\cr 0\cr 0\cr 0},
\ee

\bel{bm}
Y_{V_+}^{(l)}\pmatrix{a\cr \alpha\cr b\cr c\cr \delta\cr d} =
\pmatrix{0\cr a/2\cr -\alpha /2\cr 0\cr c/2\cr -\delta /2},\qquad
X_{V_+}^{(l)}\pmatrix{a\cr \alpha\cr b\cr c\cr \delta\cr d} =
\pmatrix{\gamma /2\cr e/2\cr \beta /2\cr 0\cr 0\cr 0},
\ee

\bel{bn}
Y_{V_-}^{(r)}\pmatrix{a\cr \alpha\cr b\cr c\cr \delta\cr d} =
\pmatrix{-\alpha /2\cr b/2\cr 0\cr -\delta /2\cr d/2\cr 0},\qquad
X_{V_-}^{(r)}\pmatrix{a\cr \alpha\cr b\cr c\cr \delta\cr d} =
\pmatrix{0\cr 0\cr 0\cr -\gamma /2\cr e/2\cr -\beta /2},
\ee

\bel{bo}
Y_{V_-}^{(l)}\pmatrix{a\cr \alpha\cr b\cr c\cr \delta\cr d} =
\pmatrix{\alpha /2\cr b/2\cr 0\cr \delta /2\cr d/2\cr 0},\qquad
X_{V_-}^{(l)}\pmatrix{a\cr \alpha\cr b\cr c\cr \delta\cr d} =
\pmatrix{0\cr 0\cr 0\cr \gamma /2\cr e/2\cr \beta /2}.
\ee
\een

\noindent
In non-coboundary cases Poisson brackets can be calculated by
applying supersymmetric version of the method described in
\cite{Sobczyk}. It is necessary to solve the cocycle equation for an
element $\Phi : G\rightarrow {\cal G}\wedge {\cal G}$

\bel{bp}
\Phi (g_1 g_2) = \Phi (g_1) + g_1 \Phi (g_2) g_1^{-1}
\ee

\noindent
satisfying "initial conditions" determined by the Lie super-bialgebra
structure in consideration.

The whole procedure has to be followed in every detail only in the case
of Lie super-bialgebra structures (i) and (iv)
from the list presented in
Chapter~\ref{sec3}. One finds

\bel{bpa}
\Phi_{(i)} (g) = cs P_+\wedge P_-
\ee

\bel{bpb}
\ba{rcl}
\Phi_{(iv)} (g) &=& -2a {\rm e}^s P_+\wedge H - a {\rm e}^s D_+\wedge D_+ 
- 2b {\rm e}^{-s} P_-\wedge H \\
&& +2ab P_-\wedge P_+ + b {\rm e}^{-s} D_-\wedge D_- + \xi {\rm e}^{s\over 2}
H\wedge D_+ \\
&&  -a\xi {\rm e}^{3s\over 2} P_+\wedge D_+
+\xi b {\rm e}^{-{s\over 2}}P_-\wedge D_+
+\eta{\rm e}^{-{s\over 2}}H\wedge D_- \\
&& -a\eta {\rm e}^{s\over 2}P_+\wedge D_-
+\eta b{\rm e}^{-{3s\over 2}} P_-\wedge D_- - {1\over 2} \xi\eta D_+\wedge D_-\,.
\ea
\ee

\noindent
For $\Phi (g) = \Phi^{jk} (g)  g_j\wedge g_k$
Poisson-Lie brackets are

\bel{bqa}
\{\phi , \psi \} =
\left( X_j^{(r)}\phi\right) \Phi^{jk} \left( X_k^{(l)}\psi\right) .
\ee

\noindent
In all the remaining cases one can use a fact that
for the special value of parameter $c$: $c=0$ they become
coboundaries.
The parameter $c$ is present only in $\delta (H)$
as $cP_+\wedge P_-$. It is clear that it
will appear in $\Phi$ as
$csP_+\wedge P_-$. It is thus possible to calculate Poisson-Lie
brackets using classical $r$-matrices given in \r{19}
and to add at the very end
of computations
the extra term in $\{a, b\}$. The complete set of
relations making $OSp(1,2)$ and super $E(2)$ Poisson-Lie supergroups
are given below in two tables. In the case of $OSp(1,2)$ 
in order to obtain shorter formulae all the Poisson
brackets have been multiplied by the factor $2$.

\def\a{\alpha}
\def\b{\beta}
\def\c{\gamma}
\def\d{\delta}
\def\k{\xi}
\def\n{\eta}
\def\cc#1{\multicolumn{1}{c|}{#1}}

$$
\begin{array}{|c|r|r|r|}
\multicolumn{4}{l}{\mbox{Table 1: Poisson Lie structures for 
the group $OSp(1|2)$}.}\\[2mm]
\hline
 &  \cc{1}  &  \cc{2} &  \cc{3} \\
 \hline
 \{a,b\} & a^2+\a\d-1 & a^2-1 &a^2+b^2-1\\  
 \{a,c\}  & -c^2          & -c^2    & 1-c^2-a^2\\
 \{a,d\} & c(a-d)         & c(a-d)  &(c-b)(a-d)\\
 \{b,c\} & -c(a+d)        & -c(a+d) &-(b+c)(a+d)\\
 \{b,d\} & 1-d^2-\a\d   & 1-d^2 &1-b^2-d^2)\\
 \{c,d\}  & c^2           & c^2     & c^2+d^2-1\\
 \{a,\a\} & c\a - a\d   &           & b\a\\
 \{b,\a\} & d\a-b\d     & -a\a    & -a\a\\
 \{c,\a\} & c\d           & c\d     & c\d+a\a+b\d\\
 \{d,\a\} & d\d           & (d-a)\d & d\d+b\a-a\d\\
 \{a,\d\} & -c\d          & -c\d    & -c\d-a\a+d\a\\
 \{b,\d\} & -d\d            & -d\d-c\a &-(d\d+b\a+c\a)\\
 \{c,\d\} &                 &           & d\d\\
 \{d,\d\}  &                 & -c\d   & -c\d\\
 \{\a,\a\} & 2\d\a           & 1-a^2 &1-a^2-b^2\\
 \{\a,\d\} &                & -ac     & -ac-bd\\
 \{\d,\d\} &                & -c^2    & 1-c^2-d^2\\
 \hline
\end{array}
$$

$$
\begin{array}{|c|r|r|r|r|r|r|}
\multicolumn{7}{l}{\mbox{Table 2: Poisson Lie structures for 
the group super-$E(2)$}.}\\[2mm]
\hline
          &  \cc{(i)} &  \cc{(ii)}   &   \cc{(iii)}		& \cc{(iv)}    &
           \cc{(v)} &            \cc{(vi)}  \\
 \hline
 \{a,b\}   &cs & -b+cs &a-b+cs         &-2ab    &-b-e^s   & a-b+cs       \\
 \{a,e^s\} &    & 1-e^s	&1-e^s 		&-2ae^s  &1-e^s    & 1-e^s        \\  
 \{b,e^s\} &    & 	&e^s-e^{2s} 	&-2be^s  &         & e^s - e^{2s}\\  
 \{a,\k\}  &    & \k/2	& \k/2		&        &-\k e^{-s}/2 &-\k e^{-s}/2     \\  
 \{a,\n\}  &    & -\n/2	&\k-(\n/2)	&-a\n    &-\n/2    &  -\n /2       \\  
 \{b,\k\}  &    & 	&\n-(\k/2) 	&\k b    &         & -\k/2	\\  
 \{b,\n\}  & 	& 	&\n/2 		&        &         & -\n e^s/2  \\  
 \{e^s,\k\}& 	& 	& 		& \k e^s &         & 		\\  
 \{e^s,\n\}& 	& 	& 		& \n e^s &         & 		\\  
 \{\k,\k\} & 	& 	& 		&-2a     &e^{-s}-1 & e^{-s}-1	\\  
 \{\k,\n\} & 	& 	& 		&-\k\n/2 &         & 		\\  
 \{\n,\n\} & 	& 	& 		& 2b     &         & e^s-1	\\  
 \hline
\end{array}
$$

\section{Discussion and final remarks}\label{sec6}

In five cases of Poisson-Lie brackets on the super-$E(2)$ group it is
straightforward to obtain quantum group structures. In fact, RHS's of
Poisson brackets are rather simple and contain no products of functions
which could lead to ordering ambiguities upon quantization. Deformation
of super-$E(2)$ introduced in \cite{Kontrakcja}
corresponds in our classification to the structure (iv).
By applying
duality techniques in cases (i), (ii), (iii), (v) and (vi)
it must be possible to obtain quantum deformations
of super $e(2)$ algebra. 
In the case of $OSp(1,2)$ a quantization of the first
structure appeared very recently in \cite{CeKu}. The second structure is
the standard one introduced in \cite{Kulish}. A quantum version of the third
structure is missing.

\subsection*{Acknowledgment}

The authors would like to thank Prof.\ J.\ Lukierski for valuable discussions.

\end{document}